\documentclass[a4paper,12pt]{article}
\usepackage[utf8]{inputenc}
\usepackage{amsmath}
\usepackage{amssymb}
\usepackage{mathtools}
\usepackage{fancyhdr}
\usepackage{siunitx}
\usepackage{mleftright}
\usepackage{booktabs}
\usepackage{tabularx}
\usepackage{caption}
\usepackage{aligned-overset}
\usepackage{subcaption}
\usepackage[intlimits]{esint}
\usepackage{stackengine}
\usepackage[explicit]{titlesec}
\usepackage{tocloft}
\usepackage{lipsum}
\usepackage{authblk}
\usepackage[textwidth=15cm,textheight=22cm]{geometry}

\renewcommand{\headrulewidth}{0.4pt}

\DeclareMathAlphabet{\mathbbold}{U}{bbold}{m}{n}

\usepackage[hyphens,spaces,obeyspaces]{url}
\usepackage[pdftitle={Analogs to Ramanujan's Master Theorem}, pdfauthor={Julius Lehmann},colorlinks=true, citecolor=blue, urlcolor=black, linkcolor=blue]{hyperref}
\usepackage[style=numeric]{biblatex}
\newcommand{\dr}[1]{\mathrm{#1}}
\newcommand{\dd}[1]{\dr{d}#1\!\mathop{}}
\newcommand{\der}[1]{\frac{\mathrm{d}}{\mathrm{d}#1}}

\newcommand{\DerN}[3]{\frac{\mathrm{d}^{#3}#1}{\mathrm{d}{#2}^{#3}}}

\newcommand{\nexp}[1]{\dr{e}^{#1}}

\newcommand{\x}{\dr{i}}

\DeclarePairedDelimiter\klam{(}{)}
\DeclarePairedDelimiter\reklam{[}{]}

\DeclarePairedDelimiterX\braket[2]{\langle}{\rangle}{#1 \delimsize\vert #2}

\title{Analogs to Ramanujan's Master Theorem and Operational Methods}
\author{Julius Lehmann\footnote{e-mail: \texttt{julius.lehmann(at)tum.de}}}
\affil{\textit{Physics of Complex Biosystems, Technical University of Munich, 85748 Garching, Germany}}
\date{\today}

\normalsize

\begin{document}

\maketitle

\begin{abstract}
   In this paper, we utilize operational methods to obtain closed-form solutions for certain classes of integrals in the spirit of Ramanujan's Master Theorem and provide several analogs to it. Although the use of operational calculus makes the proofs formal in nature, they can still yield interesting and correct results and may stimulate further rigorous investigations in the future.
\end{abstract}

\section{Ramanujan's Master Theorem}
Since the advent of calculus, there has been particular interest in finding closed-form solutions to general classes of integrals. This problem has been further stimulated by the fact that derivatives can be computed in an algorithmic manner, while even simple looking integrals can require a considerable amount of work and ingenuity to evaluate. This considered, it is remarkable that it is possible to find a rather simple closed-form solution to a class of integrals, namely Mellin transforms, using \emph{Ramanujan's master theorem} (RMT), first described by Srinivasa Ramanujan (1887-1920) at the beginning of the 20${}^\text{th}$ century. Before giving the expression of his master theorem, we note that a similar result in the spirit of the RMT already appeared in the 19${}^\text{th}$ century in a paper by J.\,W.\,L. Glaisher~\cite{glaisher1874vii,amdeberhan2012ramanujan}. By integrating term-by-term an identity related to the Newton interpolation formula of a function, he finds~\cite{FootnoteNotation}
\begin{equation}\label{eq:Glaisher}
    \int_0^\infty\dd{x}\sum_{n=0}^\infty(-1)^n\phi(n)x^{2n}=\int_0^\infty\dd{x}(\phi(0)-\phi(1)x^2+\phi(2)x^4+\dots)=\frac{\pi}{2}\phi\klam[\Big]{-\frac{1}{2}}
\end{equation}
to which he remarks, ``[...] of course, $\phi(n)$ being only defined for $n$ a positive integer, $\phi\klam[\big]{-\frac{1}{2}}$ is without meaning. But in cases where $\phi(n)$ involves factorials, there is a strong presumption [...] that the formula will give correct result if the continuity of the terms is preserved by the substitution of gamma functions for the factorials. This I have found to be true in every case to which I have applied \eqref{eq:Glaisher}.''~\cite{glaisher1874vii} It does not appear that much further work was done along these lines apart from a paper by J. O'Kinealy~\cite{o1874new} which was well  received by J. Glaisher, until almost 40 years later, Ramanujan gives the following theorem in his First Quarterly Report~\cite{berndt1983quarterly,berndt2012ramanujan,hardy1999ramanujan}: 
\\\\
\noindent\textbf{Ramanujan's Master Theorem.} \textit{Suppose that in some neighborhood around $x=0$ we can expand the function $f(x)$ as a power series of the form
\begin{equation*}
    f(x)=\sum_{n=0}^\infty\frac{(-1)^n\phi(n)}{n!}x^n,
\end{equation*}
with $\phi(0)\neq0$, then under certain growth conditions on $\phi$ the following relation holds:
\begin{equation*}\label{eq:RMT}\tag{RMT}
    \int_0^\infty\dd{x}x^{s-1}f(x)=\Gamma(s)\phi(-s),\quad0<s
\end{equation*}
where $\phi(-s)$ is interpreted as a ``natural'' analytic continuation to the real numbers of the discrete sequence $\{\phi(n)\}_{n\in\mathbb{N}_0}$.}
\\\\
This theorem was widely used by him as a tool in computing definite integrals and infinite series. Indeed, as G.\,H. Hardy puts it in~\cite{hardy1999ramanujan}, he ``was particularly fond of them, and used them as one of his commonest tools.'' Since many examples of the application of \ref{eq:RMT} are available in the literature and by Ramanujan himself, we will refer to Refs.~\cite{amdeberhan2012ramanujan,berndt2012ramanujan,atale2021certain} for a collection of nice examples. It is possible to rigorously prove the RMT using complex analysis and the inversion formula of the Mellin transform under certain, as G.\,H. Hardy describes them, ``natural'' conditions on $\phi$~\cite{hardy1999ramanujan}. However, similar to how Ramanujan himself proceeded, it can be proven quite easily using formal operational methods. This will be shown in the next paragraph.

To formally (albeit non-rigorously) proof the RMT, we employ the formalism of \emph{operational calculus}~\cite{OperationalCalculus}, which is a technique for manipulating operators independently of the function to which they are applied by treating them algebraically. Although its history goes back to the foundations of calculus with Gottfried Wilhelm Leibniz, the approach was mainly developed during the 19${}^\text{th}$ century by British and Irish mathematicians, including Sylvester, Boole, Glaisher, Crofton, and Blizard~\cite{bell1938history,cartier2000mathemagics}, and fully developed by Oliver Heaviside in the 1890s~\cite{nahin2002oliver}. The proof is as follows: by using the forward shift operator $\hat{\lambda}$ with the defining property $\hat{\lambda}^n\phi(\mathbbold{0})=\phi(n)$, we can rewrite the power series expansion of $f(x)$ as
\begin{equation}
    f(x)=\sum_{n=0}^\infty\frac{(-1)^n\phi(n)}{n!}x^n=\sum_{n=0}^\infty\frac{(-1)^n}{n!}(x\hat{\lambda})^n\phi(\mathbbold{0})=\nexp{-x\hat{\lambda}}\phi(\mathbbold{0}),
\end{equation}
where we used $\mathbbold{0}$ instead of just $0$ to further emphasize that $\phi$ is being acted upon. Inserting this into the integral yields
\begin{equation}
    \int_0^\infty\dd{x}x^{s-1}\nexp{-x\hat{\lambda}}\phi(\mathbbold{0})=\Gamma(s)\hat{\lambda}^{-s}\phi(\mathbbold{0})=\Gamma(s)\phi(-s),
\end{equation}
where we used the standard representation of the gamma function and, as a key step, formally treated the operator $\hat{\lambda}$ as a number. It is this algebraic manipulation, along with the implicit assumption
\begin{equation}
    \int_0^\infty\dd{x}x^{s-1}\klam[\big]{\nexp{-x\hat{\lambda}}\phi(\mathbbold{0})}=\klam[\Bigg]{\int_0^\infty\dd{x}x^{s-1}\nexp{-x\hat{\lambda}}}\phi(\mathbbold{0}),
\end{equation}
which makes the proof formal: obviously, the coefficients $\phi(n)$ must satisfy certain conditions for the result to be valid.

By setting $\phi(n)=\psi(n)\Gamma(n+1)$ in \ref{eq:RMT} we get the slightly different integral identity
\begin{equation}
    \int_0^\infty\dd{x}x^{s-1}\sum_{n=0}^\infty(-1)^n\psi(n)x^n=\frac{\pi}{\sin(\pi s)}\psi(-s),
\end{equation}
which is given by Ramanujan next to the RMT in his quarterly report~\cite{berndt1983quarterly}. Alternatively, we could have derived it again using the shift operator and noting that
\begin{equation}
    \sum_{n=0}^\infty(-1)^n\psi(n)x^n=\sum_{n=0}^\infty(-1)^n(x\hat{\lambda})^n\psi(\mathbbold{0})=\frac{1}{1+x\hat{\lambda}}\psi(\mathbbold{0})
\end{equation}
and integrate the expression. This also immediately shows that Glaisher's relation readily follows from the ansatz
\begin{equation}
    \sum_{n=0}^\infty(-1)^n\psi(n)x^{2n}=\frac{1}{1+x^2\hat{\lambda}}\psi(\mathbbold{0}).
\end{equation}

Employing the shift operator $\hat{\lambda}$, we can derive a multitude of other analogs to Ramanujan's master theorem (for related works, see~\cite{atale2021certain,bradshaw2023operational}). To begin, consider the function
\begin{equation}\label{eq:Fun1}
    f(x)=\sum_{n=0}^\infty\frac{(-1)^n\phi(n)}{n!}x^{qn+p}=x^p\nexp{-x^q\hat{\lambda}}\phi(\mathbbold{0}),
\end{equation}
with the last equality following as usual. The resulting Ramanujan-type relation reads
\begin{equation}\label{eq:RMTGen1}
    \int_0^\infty\dd{x}x^{s-1}f(x)=\frac{1}{q}\Gamma\klam[\Big]{\frac{p+s}{q}}\phi\klam[\Big]{-\frac{p+s}{q}}.
\end{equation}
Obviously, setting $\psi(n)=\phi(n)\Gamma(n+1)$ gives the corresponding RMT for non-exponential power series. A similar result is given in Refs.~\cite{atale2021certain,gorska2012ramanujan}. A nice application of Eq.~\eqref{eq:RMTGen1} is the following integral, where $J_\alpha(x)$ is the Bessel function of the first kind of order $\alpha$:
\begin{equation}
    \int_0^\infty\dd{x}x^{s-1}J_\alpha(x^m).
\end{equation}
Recalling the series definition of the Bessel function~\cite{abramowitz1968handbook},
\begin{equation}
    J_\alpha(x)=\sum_{n=0}^\infty\frac{(-1)^n}{n!\Gamma(n+\alpha+1)}\klam[\Big]{\frac{x}{2}}^{2n+\alpha},
\end{equation}
$J_\alpha(x^m)$ is represented by the choices $\phi(n)=2^{-2n-\alpha}/\Gamma(n+\alpha+1)$, $q=2m$, and $p=\alpha m$. This gives us the expression for the integral of
\begin{equation}
    \int_0^\infty\dd{x}x^{s-1}J_\alpha(x^m)=2^{\frac{s}{m}-1}\frac{\Gamma(\frac{\alpha}{2}+\frac{s}{2m})}{m\Gamma(1+\frac{\alpha}{2}-\frac{s}{2m})}.
\end{equation}
This represents a generalization of entry 6.561.14 of~\cite{gradshteyn2014table} to higher powers of the argument of $J_\alpha$. Another interesting analog which was first derived by Atale et al.~\cite{atale2021certain} involves the Riemann zeta function defined by
\begin{equation}
    \zeta(s)=\sum_{n=1}^\infty\frac{1}{n^s}.
\end{equation}
To derive it, take the \ref{eq:RMT} and replace $x\to nx$:
\begin{equation}
    \int_0^\infty\dd{x}x^{s-1}n^sf(nx)=\Gamma(s)\phi(-s).
\end{equation}
Dividing both sides by $n^s$ and then summing over all $n$ gives
\begin{equation}
    \int_0^\infty\dd{x}x^{s-1}Z(x)=\zeta(s)\Gamma(s)\phi(-s),
\end{equation}
where we defined
\begin{equation}
    Z(x)\equiv\sum_{n=1}^\infty f(nx).
\end{equation}
An immediate consequence of this expression is the standard integral definition of the Riemann zeta function~\cite{apostol2013introduction}:
\begin{equation}
    \int_0^\infty\dd{x}\frac{x^{s-1}}{\nexp{x}-1}=\zeta(s)\Gamma(s)
\end{equation}
which follows from $Z(x)=\sum_{n=1}^\infty\nexp{-nx}$ implying $\phi(n)=1$.

The third analog is a direct consequence of the RMT together with the relation $\ln(x)=\der{s}x^s|_{s=0}$. This allows us to evaluate general logarithmic integrals of the form
\begin{equation}
    \int_0^\infty\dd{x}x^{s-1}\ln^m(x)f(x),
\end{equation}
where $m$ is a positive integer and the function has a series expansion of type Eq.~\eqref{eq:Fun1}. By taking successive derivatives with respect to $s$, we obtain the logarithmic analog of Ramanujan's master theorem:
\begin{equation}\label{eq:RMTLog}
    \int_0^\infty\dd{x}x^{s-1}\ln^m(x)f(x)=\frac{1}{q}\DerN{}{s}{m}\reklam[\bigg]{\Gamma\klam[\Big]{\frac{p+s}{q}}\phi\klam[\Big]{-\frac{p+s}{q}}},
\end{equation}
of which we have $m=1$ as a useful special case:
\begin{equation}
    \int_0^\infty\dd{x}x^{s-1}\ln(x)f(x)=\frac{1}{q^2}\Gamma\klam[\Big]{\frac{p+s}{q}}\reklam[\bigg]{\psi_0\klam[\Big]{\frac{p+s}{q}}\phi\klam[\Big]{-\frac{p+s}{q}}-\phi'\klam[\Big]{-\frac{p+s}{q}}}.
\end{equation}
$\psi_0(x)\equiv\der{x}\ln(\Gamma(x))$ in the above equation is the digamma function~\cite{abramowitz1968handbook}. A simple application of Eq.~\eqref{eq:RMTLog} is the well-known expression
\begin{equation}
    \int_0^\infty\dd{x}x^{s-1}\ln^2(x)\nexp{-x}=\DerN{}{s}{2}\Gamma(s)=\Gamma(s)\bigl(\psi_0^2(s)+\psi_1(s)\bigr),
\end{equation}
where $\psi_1(x)$ is the trigamma function, the derivative of the digamma function~\cite{abramowitz1968handbook}. A more interesting example, however, is the following integral:
\begin{equation}
    \int_0^\infty\dd{x}\frac{x^{s-1}\ln^m(x)}{(1+x^\alpha)^\beta}.
\end{equation}
Notice that our function can be written as
\begin{equation}
    \frac{1}{(1+x^\alpha)^\beta}=\sum_{n=0}^\infty\frac{(-1)^n}{n!}\frac{\Gamma(\beta+n)}{\Gamma(\beta)}x^{\alpha n},
\end{equation}
we can deduce $\phi(n)=\Gamma(\beta+n)/\Gamma(\beta)$, $q=\alpha$, and $p=0$. Therefore, using the logarithmic analog to the RMT, we find
\begin{equation}
    \int_0^\infty\dd{x}\frac{x^{s-1}\ln^m(x)}{(1+x^\alpha)^\beta}=\frac{1}{\alpha\Gamma(\beta)}\DerN{}{s}{m}\reklam[\bigg]{\Gamma\klam[\Big]{\frac{s}{\alpha}}\Gamma\klam[\Big]{\beta-\frac{s}{\alpha}}},
\end{equation}
so for example, $m=1$:
\begin{equation}
    \int_0^\infty\dd{x}\frac{x^{s-1}\ln(x)}{(1+x^\alpha)^\beta}=\frac{\Gamma\klam{\frac{s}{\alpha}}\Gamma\klam{\beta- \frac{s}{\alpha}}}{\alpha\Gamma(\beta)}\reklam[\Big]{\psi_0\klam[\Big]{\frac{s}{\alpha}}-\psi_0\klam[\Big]{\beta-\frac{s}{\alpha}}},
\end{equation}
or taking $m=2$, $s=\frac{1}{2}$, and $a=b=1$:
\begin{equation}
    \int_0^\infty\dd{x}\frac{\ln^2(x)}{\sqrt{x}(1+x)}=\pi^3.
\end{equation}

One takeaway from the formal proof of RMT that allows us to generalize the result fairly easily is that whenever the series expansion of a function is given by
\begin{equation}
    f(x)=\sum_{n=0}^\infty\phi(n)\gamma(n)x^n,
\end{equation}
where $\gamma(n)$ are the expansion coefficients of another function $g(x)$, so that we can write $f(x)$ in terms of the shift operator as
\begin{equation}\label{eq:gRMT}
    f(x)=\sum_{n=0}^\infty\psi(n)(x\hat{\lambda})^n\phi(\mathbbold{0})=g(x\hat{\lambda})\phi(\mathbbold{0}),
\end{equation}
the generalized Ramanujan Master Theorem (gRMT) for an integral kernel $K(x)$ is given by
\begin{equation}
    \int_0^\infty\dd{x}K(x)f(x)=\int_0^\infty\dd{x}K(x)g(x\hat{\lambda})\phi(\mathbbold{0}).
\end{equation}
We can formulate this statement in a different, operational way: by representing the RHS as an operator, we can write
\begin{equation}
    \int_0^\infty\dd{x}K(x)f(x)=\hat{\mathcal{R}}\phi(\mathbbold{0})
\end{equation}
where now all the relevant information for evaluating the integral is contained in the operator $\hat{\mathcal{R}}\equiv\int_0^\infty\dd{x}K(x)g(x\hat{\lambda})$, which acts on the coefficients of the expansion of $f(x)$. This is in accordance with similar results given in~\cite{bradshaw2023operational}. The \ref{eq:RMT} is given as a special case by the choices $K(x)=x^{s-1}$ and $g(x)=\nexp{-x}$ or $g(x)=\frac{1}{1+x}$, respectively. Now let us consider another choice of $g(x)$. In particular, take the trigonometric functions $g_\dr{c}(x)=\cos(x)$ and $g_\dr{s}(x)=\sin(x)$. To find the expansion coefficients $\gamma_\dr{c}(n)$, consider the power series of $\cos(x)$ and rewrite it accordingly:
\begin{equation}
    \cos(x)=\sum_{n=0}^\infty\frac{(-1)^n}{(2n)!}x^{2n}=\sum_{n=0}^\infty\frac{\cos(\frac{\pi}{2}n)}{n!}x^n,
\end{equation}
since $\cos(\frac{\pi}{2}n)=(-1)^n$ for $n$ even and $0$ for $n$ odd. Therefore, we find $\gamma_\dr{c}(n)=\cos(\frac{\pi}{2}n)/n!$, and similarly, $\gamma_\dr{s}(n)=\sin(\frac{\pi}{2}n)/n!$. The resulting trigonometric RMTs that follow from Eq.~\eqref{eq:gRMT} are therefore given by
\begin{equation}
    \begin{aligned}
        \int_0^\infty\dd{x}x^{s-1}f_\dr{c}(x)&=\int_0^\infty\dd{x}x^{s-1}\sum_{n=0}^\infty\frac{(-1)^n\phi(2n)}{(2n)!}x^{2n}=\Gamma(s)\cos\klam[\Big]{\frac{\pi}{2}s}\phi(-s)\\
    \int_0^\infty\dd{x}x^{s-1}f_\dr{s}(x)&=\int_0^\infty\dd{x}x^{s-1}\sum_{n=0}^\infty\frac{(-1)^n\phi(2n+1)}{(2n+1)!}x^{2n+1}=\Gamma(s)\sin\klam[\Big]{\frac{\pi}{2}s}\phi(-s).
    \end{aligned}
\end{equation}
An interesting observation taken from the gRMT for $\cos$- or $\sin$-type series is the following: by using the main assumption of operational calculus, i.\,e. treating operators algebraically, we can identify the associated $\hat{\mathcal{R}}$ as simply being the Fourier cosine/sine transform of $K(x)$ in $\hat{\lambda}$~\cite{coleman2016introduction}, namely:
\begin{equation}
    \begin{aligned}
        \int_0^\infty\dd{x}K(x)\cos(x\hat{\lambda})\phi(\mathbbold{0})&=\frac{\pi}{2}\mathcal{F}_\dr{c}\{K(x)\}(\hat{\lambda})\phi(\mathbbold{0})\\
        \int_0^\infty\dd{x}K(x)\sin(x\hat{\lambda})\phi(\mathbbold{0})&=\frac{\pi}{2}\mathcal{F}_\dr{s}\{K(x)\}(\hat{\lambda})\phi(\mathbbold{0}).
    \end{aligned}
\end{equation}
Since there are many transformation pairs available in the literature, this opens the door to a plethora of Ramanujan-type master theorems involving $\cos$- or $\sin$-type functions. Consider the case $K(x)=\nexp{-x^2}$, then we have
\begin{equation}
    \int_0^\infty\dd{x}\nexp{-x^2}f_\dr{c}(x)=\frac{\pi}{2}\mathcal{F}_\dr{c}\{\nexp{-x^2}\}(\hat{\lambda})\phi(\mathbbold{0})=\frac{\sqrt{\pi}}{2}\nexp{-\frac{\hat{\lambda}^2}{4}}\phi(\mathbbold{0}).
\end{equation}
By interpreting the RHS as a formal power series in $\hat{\lambda}$, we arrive at the following power series solution:
\begin{equation}
    \int_0^\infty\dd{x}\nexp{-x^2}f_\dr{c}(x)=\frac{\sqrt{\pi}}{2}\sum_{n=0}^\infty\frac{(-1)^n\phi(2n)}{4^nn!},
\end{equation}
and analogously for a $\sin$-type series.

In the penultimate part of this section, we apply the general Ramanujan's master theorem to integral transformations. The original form of the RMT deals with Mellin transforms, other very important transforms are the already mentioned Fourier cosine and sine transforms, the Laplace and Hankel transform~\cite{bracewell1966fourier,piessens2000hankel}. The former are generated by choosing the kernels $K(x)=\cos(sx)$ and $K(x)=\sin(sx)$, while the latter have kernels $K(x)=\nexp{-sx}$ and $K(x)=J_\alpha(sx)x$. We start with the Fourier-type transforms:
\begin{equation}
    \mathcal{F}_\dr{c}\{f(x)\}(s)=\int_0^\infty\dd{x}\cos(sx)f(x),\qquad\mathcal{F}_\dr{s}\{f(x)\}(s)=\int_0^\infty\dd{x}\sin(sx)f(x).
\end{equation}
By assuming an exponential-type power series expansion for $f(x)$, we find the solutions in terms of $\hat{\lambda}$ and the expansion coefficients $\phi(\mathbbold{0})$:
\begin{equation}
    \begin{aligned}
        \int_0^\infty\dd{x}\cos(sx)f(x)&=\frac{\hat{\lambda}}{s^2+\hat{\lambda}^2}\phi(\mathbbold{0})\\
        \int_0^\infty\dd{x}\sin(sx)f(x)&=\frac{s}{s^2+\hat{\lambda}^2}\phi(\mathbbold{0}).
    \end{aligned}
\end{equation}
What remains is to define the action of the operators acting on $\phi$. We again interpret the fractions as formal power series, but we have two ways to expand them: either as series around the origin or as asymptotic series at infinity. Both approaches yield series in powers of $s$, one being a small $s\approx0$ expansion, the other a large $s\to\infty$ expansion. Taking this into account, we find
\begin{equation}
    \begin{aligned}
        \int_0^\infty\dd{x}\cos(sx)f(x)&=\sum_{n=0}^\infty(-1)^n\phi(-2n-1)s^{2n}\\
        \int_0^\infty\dd{x}\cos(sx)f(x)&\sim\sum_{n=0}^\infty(-1)^n\frac{\phi(2n+1)}{s^{2n}}.
    \end{aligned}
\end{equation}
and
\begin{equation}
    \begin{aligned}
        \int_0^\infty\dd{x}\sin(sx)f(x)&=\sum_{n=0}^\infty(-1)^n\phi(-2n)s^{2n+1}\\
        \int_0^\infty\dd{x}\sin(sx)f(x)&\sim\sum_{n=0}^\infty(-1)^n\frac{\phi(2n)}{s^{2n+1}}.
    \end{aligned}
\end{equation}
Without going into too much detail, following the same procedure, we find the analogous expressions for the Laplace transform
\begin{equation}
    \mathcal{L}\{f(x)\}(s)=\int_0^\infty\dd{x}\nexp{-sx}f(x)
\end{equation}
given by the formulas
\begin{equation}
    \begin{aligned}
        \int_0^\infty\dd{x}\nexp{-sx}f(x)&=\sum_{n=0}^\infty(-1)^n\phi(-n-1)s^n\\
        \int_0^\infty\dd{x}\nexp{-sx}f(x)&\sim\sum_{n=0}^\infty(-1)^n\frac{\phi(n)}{s^{n+1}}.
    \end{aligned}
\end{equation}
Note that the second expression above is the well-known asymptotic expansion of the Laplace transform found by term-wise integration of the power series expansion of $f(x)$. The Hankel transform frequently arises in in the evaluation of Fourier transforms of radially symmetric functions. It is defined by
\begin{equation}
    H_\alpha\{f(x)\}(s)=\int_0^\infty\dd{x}xJ_\alpha(sx)f(x)
\end{equation}
and gives, under the usual assumptions, rise to the following associated operator in terms of $\hat{\lambda}$:
\begin{equation}
    H_\alpha\{f(x)\}(s)=\klam[\bigg]{\frac{s}{\hat{\lambda}+\sqrt{s^2+\hat{\lambda}^2}}}^\alpha\frac{\hat{\lambda}+\alpha\sqrt{s^2+\hat{\lambda}^2}}{(s^2+\hat{\lambda}^2)^{\frac{3}{2}}}\phi(\mathbbold{0}).
\end{equation}
Since a general expression for the expansion of the RHS is very cumbersome, we will focus on the special case $\alpha=0$. There we have
\begin{equation}
    H_\alpha\{f(x)\}(s)=\frac{\hat{\lambda}}{(s^2+\hat{\lambda}^2)^{\frac{3}{2}}}\phi(\mathbbold{0}),
\end{equation}
which leads to the following formulas:
\begin{equation}
    \begin{aligned}
        \int_0^\infty\dd{x}xJ_0(sx)f(x)&=\sum_{n=0}^\infty\binom{-\frac{3}{2}}{n}\phi(-2n)s^{2n}\\\int_0^\infty\dd{x}xJ_0(sx)f(x)&\sim\sum_{n=0}^\infty\binom{-\frac{3}{2}}{n}\frac{\phi(2n+1)}{s^{2n+3}}.
    \end{aligned}
\end{equation}
The last expression was also found and rigorously derived using the Hardy-Ramanujan master theorem~\cite{hardy1999ramanujan} in Ref.~\cite{kisselev2021exact}.

Before we enter the next section, let us briefly address another relation given by Ramanujan in his Second Quarterly Report and a double-integral analog to the RMT. Alongside the other relations, we find in \cite{berndt1983quarterly} the statement (adapted for the purpose of the paper)
\begin{equation}\label{eq:FG}
    \int_0^\infty\dd{x}\klam*{\sum_{n=0}^\infty\frac{(-1)^n\phi(n)}{n!}x^n}\klam*{\sum_{n=0}^\infty\frac{(-1)^n\psi(n)}{n!}x^n}=\sum_{n=0}^\infty(-1)^n\phi(n)\psi(-n-1),
\end{equation}
which can also be given with the roles of $\phi$ and $\psi$ in the RHS reversed. This relation follows in a straightforward manner from writing the sums in their respective shift operator form and then expanding the resulting function of operators. To derive the next relation involving a double integral, we start with the well-known representation of the 0${}^\text{th}$ order modified Bessel function~\cite{abramowitz1968handbook}:
\begin{equation}
    K_0(s\xi)=\int_0^\infty\dd{x}\frac{\cos(sx)}{\sqrt{x^2+\xi^2}}.
\end{equation}
Given a function $f(x)$ with a $\cos$-type series expansion
\begin{equation}
    f(x)=\sum_{n=0}^\infty\frac{(-1)^n\phi(2n)}{(2n)!}x^{2n},
\end{equation}
we can utilize this integral to arrive at the operational statement
\begin{equation}
    \int_0^\infty\dd{x}\frac{f(x)}{\sqrt{x^2+\xi^2}}=K_0(\xi\hat{\lambda})\phi(\mathbbold{0}).
\end{equation}
Multiplying both sides by $\xi^{s-1}$ and then integrating the equation yields the expression
\begin{equation*}
    \int_0^\infty\dd{\xi}\int_0^\infty\dd{x}\frac{\xi^{s-1}f(x)}{\sqrt{x^2+\xi^2}}=\Gamma\klam[\Big]{\frac{s}{2}}^2\frac{\phi(-s)}{2^{2-s}},
\end{equation*}
which generalizes Ramanujan's original master theorem to double integrals.

\section{Ramanujan's Master Theorem for Sums}
In this brief section we want to present one of Ramanujan's results on discrete analogs of his famous master theorem~\ref{eq:RMT}. Ramanujan's interest not only covered obtaining formulas for evaluating general classes of integrals, but also extended to finding similar results for general classes of sums. In chapter 4 of his first notebook~\cite{berndt2012ramanujan}, which deals with iterations of the exponential function and Bell polynomials, he states the following three relations:
\begin{subequations}
    \begin{gather}
        \sum_{n=0}^\infty\frac{(-1)^n}{2n+1}(\phi(2n+1)+\phi(-2n-1))=\frac{\pi}{2}\phi(0)\label{eq:Sum}\\
    \sum_{n=1}^\infty(-1)^{n+1}(\phi(n)+\phi(-n))=\phi(0)\label{eq:Sum1}\\
    \sum_{n=1}^\infty\frac{(-1)^{n+1}}{n}(\phi(n)-\phi(-n))=\phi'(0)\label{eq:Sum2},
    \end{gather}
\end{subequations}
which are valid under conditions elaborated in~\cite{berndt2012ramanujan}. We will illustrate Ramanujan's derivation using the first equation~\eqref{eq:Sum}. As for the master theorem, the proof is again formal in nature and therefore non-rigorous. We start by considering the relation $\arctan(x)+\arctan(1/x)=\pi/2$ which holds for $x>0$. Suppose that in some neighborhood around $x=0$ we can expand the function $f(x)$ as an exponential power series with expansion coefficient $\psi(n)$. Replacing $x$ by some power of it in the $\arctan$-relation clearly does not change its value, therefore Ramanujan writes down the following set of equations:
\begin{equation}
    \begin{aligned}
        \frac{\psi(0)}{0!}\klam[\Big]{\arctan(x^0)+\arctan\klam[\Big]{\frac{1}{x^0}}}&=\frac{\pi}{2}\frac{\psi(0)}{0!}\\
    \frac{\psi(1)}{1!}\klam[\Big]{\arctan(x^1)+\arctan\klam[\Big]{\frac{1}{x^1}}}&=\frac{\pi}{2}\frac{\psi(1)}{1!}\\
    \frac{\psi(2)}{2!}\klam[\Big]{\arctan(x^2)+\arctan\klam[\Big]{\frac{1}{x^2}}}&=\frac{\pi}{2}\frac{\psi(2)}{2!}\\&\vdotswithin{=}\\
    \frac{\psi(n)}{n!}\klam[\Big]{\arctan(x^n)+\arctan\klam[\Big]{\frac{1}{x^n}}}&=\frac{\pi}{2}\frac{\psi(n)}{n!}.
    \end{aligned}
\end{equation}
In the next step, we take the sum over all $n$ on both sides and conclude:
\begin{equation}
    \sum_{n=0}^\infty\frac{\psi(n)}{n!}\klam[\Big]{\arctan(x^n)+\arctan\klam[\Big]{\frac{1}{x^n}}}=\frac{\pi}{2}\sum_{n=0}^\infty\frac{\psi(n)}{n!},
\end{equation}
which is still true for all $x>0$. We recognize the right side as the value of $f$ at $x=1$. Furthermore, by expanding the $\arctan$ as a power series, we find
\begin{equation}
    \sum_{n=0}^\infty\frac{\psi(n)}{n!}\sum_{k=0}^\infty\frac{(-1)^k}{2k+1}((x^n)^{2k+1}+(x^n)^{-2k-1})=\frac{\pi}{2}f(1).
\end{equation}
The step that follows is what renders this proof formal, since we (and ultimately Ramanujan) do not check for the conditions that permit it: we exchange the order of the two summations,
\begin{equation}
    \sum_{k=0}^\infty\frac{(-1)^k}{2k+1}\sum_{n=0}^\infty\frac{\psi(n)}{n!}((x^{2k+1})^n+(x^{-2k-1})^n)=\frac{\pi}{2}f(1),
\end{equation}
and write the inner sum in terms of $f(x)$:
\begin{equation}
    \sum_{k=0}^\infty\frac{(-1)^k}{2k+1}(f(x^{2k+1})+f(x^{-2k-1}))=\frac{\pi}{2}f(1).
\end{equation}
To finish the formal proof, identify $\phi(k)\equiv f(x^k)$ and note that $f(1)=f(x^0)=\phi(0)$:
\begin{equation}
    \sum_{k=0}^\infty\frac{(-1)^k}{2k+1}(\phi(2k+1)+\phi(-2k-1))=\frac{\pi}{2}\phi(0).
\end{equation}
The other equations~\eqref{eq:Sum1},\eqref{eq:Sum2} follow in a similar way by using either $x/(1+x)+(1/x)/(1+1/x)=1$ or $\ln(1+x)-\ln(1+1/x)=\ln(x)$ instead of the $\arctan$ relation. To arrive at Eq.~\eqref{eq:Sum2}, one also utilizes the observation that
\begin{equation}
    \der{s}f(x^s)=x^s\ln(x)f'(x^s).
\end{equation}
In more general terms, the first two relations rely on suitable functions that satisfy the functional equation
\begin{equation}\label{eq:FunEq}
    h(x)+h\klam[\Big]{\frac{1}{x}}=C
\end{equation}
for some $C$. In this case, assuming that $h(x)$ can be expanded as a power series around $x=0$ with coefficients $\eta(n)$, we get the general relation
\begin{equation}
    \sum_{k=0}^\infty\eta(k)(\phi(k)+\phi(-k))=C\phi(0),
\end{equation}
which should hold for some suitable set of conditions on $\phi$. In fact, there are infinitely many functions that satisfy the functional equation~\eqref{eq:FunEq}. In particular, we can choose any function $g(x)$ and set $h(x)$ equal to
\begin{equation}
    h(x)=g(x)-g\klam[\Big]{\frac{1}{x}}+\frac{C}{2},
\end{equation}
therefore $h(1/x)=g(1/x)-g(x)+C/2$ and $h(x)+h(1/x)=C$~\cite{functionalEquation}.

Two nice applications of these theorems follow by considering $\phi(n)=\nexp{\x n\theta}$, from which we have by Eq.~\eqref{eq:Sum}
\begin{equation}
    \sum_{n=0}^\infty\frac{(-1)^n\cos((2n+1)\theta)}{2n+1}=\frac{\pi}{4},
\end{equation}
which is a well-known sum, and further by Eq.~\eqref{eq:Sum2}
\begin{equation}
    \sum_{n=1}^\infty\frac{(-1)^{n+1}\sin(n\theta)}{n}=\frac{\theta}{2},
\end{equation}
which is the familiar Fourier series of the piecewise linear function. A more complicated sum can be found by substituting $\phi(n)=\sin(n\theta)/n$, where $-\pi/2<\theta<\pi/2$, in Eq.~\eqref{eq:Sum}, from which we obtain another Fourier series of the piecewise linear function (and, in fact, the Fourier series of the triangular wave for $\theta\in\mathbb{R}$):
\begin{equation}
    \sum_{n=0}^\infty\frac{(-1)^n\sin((2n+1)\theta)}{(2n+1)^2}=\frac{\pi}{4}\theta.
\end{equation}

\section*{List of Ramanujan's Master Theorem Analogs}

\begin{table}[h]
    \centering
    \begin{tabularx}{\textwidth}{c|c|c}
        \toprule Integral & Function & Closed-Form Solution\\\midrule
        $\displaystyle\int_0^\infty\dd{x}x^{s-1}f(x)$ & $\displaystyle\sum_{n=0}^\infty\frac{(-1)^n\phi(n)}{n!}x^n$ & $\displaystyle\Gamma(s)\phi(-s)$\\
        $\displaystyle\int_0^\infty\dd{x}x^{s-1}f(x)$ & $\displaystyle\sum_{n=0}^\infty(-1)^n\phi(n)x^n$ & $\displaystyle\frac{\pi}{\sin(\pi s)}\phi(-s)$\\
        $\displaystyle\int_0^\infty\dd{x}x^{s-1}f(x)$ & $\displaystyle\sum_{n=0}^\infty\frac{(-1)^n\phi(n)}{n!}x^{qn+p}$ & $\displaystyle\frac{1}{q}\Gamma\klam[\Big]{\frac{p+s}{q}}\phi\klam[\Big]{-\frac{p+s}{q}}$\\
        $\displaystyle\int_0^\infty\dd{x}x^{s-1}Z(x)$ & $\displaystyle\sum_{n=1}^\infty f(nx)$ & $\displaystyle\zeta(s)\Gamma(s)\phi(-s)$\\
        $\displaystyle\int_0^\infty\dd{x}x^{s-1}\ln^m(x)f(x)$ & $\displaystyle\sum_{n=0}^\infty\frac{(-1)^n\phi(n)}{n!}x^n$ & $\displaystyle\frac{1}{q}\DerN{}{s}{m}\reklam[\bigg]{\Gamma\klam[\Big]{\frac{p+s}{q}}\phi\klam[\Big]{-\frac{p+s}{q}}}$\\\midrule
        $\displaystyle\int_0^\infty\dd{x}x^{s-1}f_\dr{c}(x)$ & $\displaystyle\sum_{n=0}^\infty\frac{(-1)^n\phi(2n)}{(2n)!}x^{2n}$ & $\displaystyle\Gamma(s)\cos\klam[\Big]{\frac{\pi}{2}s}\phi(-s)$\\
        $\displaystyle\int_0^\infty\dd{x}x^{s-1}f_\dr{s}(x)$ & $\displaystyle\sum_{n=0}^\infty\frac{(-1)^n\phi(2n+1)}{(2n+1)!}x^{2n+1}$ & $\displaystyle\Gamma(s)\sin\klam[\Big]{\frac{\pi}{2}s}\phi(-s)$\\\bottomrule
    \end{tabularx}
\end{table}

\begin{table}[t]
    \centering
    \begin{tabularx}{\textwidth}{c|c|c}
        \toprule Integral & Function & Closed-Form Solution\\\midrule
        $\displaystyle\int_0^\infty\dd{x}\cos(sx)f(x)$ & $\displaystyle\sum_{n=0}^\infty\frac{(-1)^n\phi(n)}{n!}x^n$ & ${\displaystyle\sum_{n=0}^\infty(-1)^n\phi(-2n-1)s^{2n}}\atop{\displaystyle\sim\sum_{n=0}^\infty(-1)^n\frac{\phi(2n+1)}{s^{2n}}}$\\
        $\displaystyle\int_0^\infty\dd{x}\sin(sx)f(x)$ & $\displaystyle\sum_{n=0}^\infty\frac{(-1)^n\phi(n)}{n!}x^n$ & ${\displaystyle\sum_{n=0}^\infty(-1)^n\phi(-2n)s^{2n+1}}\atop{\displaystyle\sim\sum_{n=0}^\infty(-1)^n\frac{\phi(2n)}{s^{2n+1}}}$\\\midrule
        $\displaystyle\int_0^\infty\dd{x}\nexp{-sx}f(x)$ & $\displaystyle\sum_{n=0}^\infty\frac{(-1)^n\phi(n)}{n!}x^n$ & ${\displaystyle\sum_{n=0}^\infty(-1)^n\phi(-n-1)s^n}\atop{\displaystyle\sim\sum_{n=0}^\infty(-1)^n\frac{\phi(n)}{s^{n+1}}}$\\\midrule
        $\displaystyle\int_0^\infty\dd{x}xJ_0(sx)f(x)$ & $\displaystyle\sum_{n=0}^\infty\frac{(-1)^n\phi(n)}{n!}x^n$ & ${\displaystyle\sum_{n=0}^\infty\binom{-\frac{3}{2}}{n}\phi(-2n)s^{2n}}\atop{\displaystyle\sim\sum_{n=0}^\infty\binom{-\frac{3}{2}}{n}\frac{\phi(2n+1)}{s^{2n+3}}}$\\\midrule
        $\displaystyle\int_0^\infty\dd{x}f(x)g(x)$ & see \eqref{eq:FG} & $\displaystyle\sum_{n=0}^\infty(-1)^n\phi(n)\psi(-n-1)$\\\midrule
        $\displaystyle\int_0^\infty\dd{\xi}\int_0^\infty\dd{x}\frac{\xi^{s-1}f(x)}{\sqrt{x^2+\xi^2}}$ & $\displaystyle\sum_{n=0}^\infty\frac{(-1)^n\phi(2n)}{(2n)!}x^{2n}$ & $\displaystyle\Gamma\klam[\Big]{\frac{s}{2}}^2\frac{\phi(-s)}{2^{2-s}}$\\\bottomrule
    \end{tabularx}
\end{table}

\begin{table}[t]
    \centering
    \begin{tabularx}{\textwidth}{c|c}
        \toprule Sum & \hspace{10ex}Closed-Form Solution\\\midrule
        $\displaystyle\sum_{n=0}^\infty\frac{(-1)^n}{2n+1}(\phi(2n+1)+\phi(-2n-1))$  & \hspace{7ex}$\displaystyle\frac{\pi}{2}\phi(0)$\\
        $\displaystyle\sum_{n=1}^\infty(-1)^{n+1}(\phi(n)+\phi(-n))$  & \hspace{7ex}$\displaystyle\phi(0)$\\
        $\displaystyle\sum_{n=1}^\infty\frac{(-1)^{n+1}}{n}(\phi(n)-\phi(-n))$  & \hspace{7ex}$\displaystyle\phi'(0)$\\\bottomrule
    \end{tabularx}
\end{table}

\clearpage

\lhead{}
\renewcommand{\headrulewidth}{0.0pt}
% \nocite{*}
% \setlength\bibitemsep{1.5ex}
\printbibliography

\end{document}